\documentclass{article}
\usepackage{amssymb,amsmath,amsfonts,epsfig,latexsym}

\begin{document}

\title{Asymptotic cohomological functions of toric divisors}

\author{M. Hering\thanks{University of Michigan, Department of
Mathematics, 2704 East Hall, Ann Arbor, MI 48109.  mhering@umich.edu}, 
A. K\"{u}ronya\thanks{Supported by the Leibniz program of the DFG.  Institut f\"{u}r Mathematik, Universit\"{a}t Duisburg-Essen, 
Standort Essen, Germany. alex.kueronya@uni-duisburg-essen.de}, 
S. Payne\thanks{Corresponding author.  Supported by a NSF Graduate Research 
Fellowship.  University of Michigan, Department of
Mathematics, 2704 East Hall, Ann Arbor, MI 48109. sdpayne@umich.edu} }

\maketitle

\begin{abstract}
    We study functions on the class group of a toric variety measuring
    the rates of growth of the cohomology groups of multiples of
    divisors.  We show that these functions are piecewise polynomial
    with respect to finite polyhedral chamber decompositions.  As
    applications, we express the self-intersection number of a
    $T$-Cartier divisor as a linear combination of the volumes of the
    bounded regions in the corresponding hyperplane arrangement and
    prove an asymptotic converse to Serre vanishing.
\end{abstract}

\newtheorem{theorem}{Theorem}
\newtheorem{definition}{Definition}
\newtheorem{proposition}{Proposition}
\newtheorem{lemma}{Lemma}
\newtheorem{corollary}{Corollary}

\def\vol{\mathrm{vol \, }}
\def\h{\ensuremath{\widehat{h}}}

\def\N{\ensuremath{\mathbb{N}}}
\def\Z{\ensuremath{\mathbb{Z}}}
\def\Q{\ensuremath{\mathbb{Q}}} 
\def\P{\ensuremath{\mathbb{P}}}
\def\C{\ensuremath{\mathbb{C}}}
\def\R{\ensuremath{\mathbb{R}}}
\def\O{\ensuremath{\mathcal{O}}}
\def\div{\mathrm{div}}
\def\mult{\mathrm{mult}}
\def\Pic{\ensuremath{\mathrm{Pic}}}

\def\<{\ensuremath{\langle}}
\def\>{\ensuremath{\rangle}}

\section{Introduction}

Suppose $D$ is an ample divisor on an $n$-dimensional algebraic
variety.  The sheaf cohomology of $\O(D)$ does not necessarily reflect
the positivity of $D$; $\O(D)$ may have few global sections and its
higher cohomology groups may not vanish.  However, for $m \gg 0$,
$\O(mD)$ is globally generated and all of its higher cohomology groups
vanish.  Moreover, the rate of growth of the space of global sections
of $\O(mD)$ as $m$ increases carries information on the positivity of
$D$.  Indeed, if we write $h^{0}(mD)$ for the dimension of $H^{0}(X,
\O(mD))$, then by asymptotic Riemann-Roch \cite[Example 1.2.19]{PAG},
\[
    (D^{n}) = \lim_{m} \frac {h^{0}(mD)} {m^{n}/n!}.
\]
In general, when $D$ is not necessarily ample, this limit exists and
is called the volume of $D$.  It is written $\widehat{h}^{0}(D)$ or
$\vol(D)$.  The regularity of the rate of growth of the cohomology
groups of $\O(mD)$ for $m \gg 0$ contrasts with the subtlety of the
behavior of the cohomology of $\O(D)$ itself and motivates the study
of asymptotic cohomological functions of divisors.

Lazarsfeld has shown that the volume of a Cartier divisor depends only
on its numerical equivalence class and that the volume function
extends to a continuous function on $N^{1}(X)_{\R}$ \cite[2.2.C]{PAG}. 
The volume function is polynomial on the ample cone, where it agrees
with the top intersection form.  In some special cases, including for
toric varieties, smooth projective surfaces, abelian varieties, and
generalized flag varieties, the volume function is piecewise
polynomial with respect to a locally finite polyhedral chamber
decomposition of the interior of the effective cone.  The behavior of
the volume function outside the ample cone is known to be more
complicated in general \cite{BKS}.  In this paper, we study the volume
function and its generalizations, the higher asymptotic cohomological
functions, in the toric case.

Let $X = X(\Delta)$ be a complete $n$-dimensional toric variety.  Let
$D$ be a $T$-Weil\footnote{Asymptotic cohomological functions of
Cartier divisors on general (not necessarily toric) varieties have
been studied in \cite{BKS} \cite{Kuronya} and \cite[2.2.C]{PAG}.  In
this paper, we allow Weil divisors wherever possible.  Our approach is
self-contained, and does not rely on results from the general theory.}
divisor on $X$, and $P_{D}$ the associated polytope in $M_{\R}$. 
Since $h^{0}(mD)$ is the number of lattice points in $P_{mD}$, and
since $P_{mD} = m P_{D}$ for all positive integers $m$,
$\widehat{h}^{0}(D)$ is the volume of $P_{D}$, normalized so that the
smallest lattice simplex in $M_{\R}$ has unit volume.  Oda and Park
describe, in combinatorial language, a finite polyhedral chamber
decomposition of the effective cone in the divisor class group of a
toric variety, which they call the Gelfand-Kapranov-Zelevinsky (or
GKZ) decomposition, such that the combinatorial structure of $P_{D}$
is constant as $D$ varies within each chamber \cite{OP}.  It follows
that $\vol P_{D}$ is polynomial on each of these chambers; see \cite[VIII.5
Problem 10]{Barvinok} and \cite[Example 2.8]{ELMNP}.  In particular, $\h^{0}$ extends to a
continuous, piecewise polynomial function with respect to a finite
polyhedral decomposition of the effective cone in $A_{n-1}(X)_{\R}$. 
The piecewise polynomial behavior of $\h^{0}$ also follows from
\cite[Proposition 4.13]{AIBL}, since toric varieties have ``finitely
generated linear series.''  The GKZ decomposition also arises as the
decomposition of the effective cone into ``Mori chambers'' and
``variation of GIT chambers,'' see \cite{Mori Dream Spaces}.  In an
appendix, we give a brief, self-contained account of GKZ
decompositions in the language of toric divisors.

Generalizing the volume function, we define higher asymptotic
cohomological functions of toric divisors by
\[
    \h^{i}(D) = \lim_{m} \frac {h^{i}(mD)}
    {m^{n}/n!},
\]
where $h^{i}(mD)$ is the dimension of $H^{i}(X, \O(mD))$.\footnote{The second author has studied higher asymptotic
cohomological functions of line bundles on general varieties, defined
similarly but with a $\limsup$ instead of a limit, and has shown that
these extend to continuous functions on $N^{1}(X)_{\R}$
\cite{Kuronya}.  It is not known whether the limits exist in general.}
In the toric case, it follows from local cohomology computations of
Eisenbud, Musta\c{t}\v{a}, and Stillman \cite{EMS} that there is a
decomposition of $M_{\R}$ into finitely many polyhedral regions such
that the dimensions of the graded pieces $H^{i} (X, \O(D))_{u}$ are
constant for lattice points $u$ in each region.  The regions are
indexed by collectins of rays $I \subset \Delta(1)$, and for $D = \sum
d_{\rho} D_{\rho}$, they are given by
\[
   P_{D,I} := \{ u \in M_{\R} : \< u, v_{\rho} \> \geq -d_{\rho} \mbox{
   if and only if } \rho \in I \}.
\]
In particular, the regions for $mD$ are the $m$-fold dilations of the
regions for $D$.  In Section 3, we deduce from this that the limit in
the definition of $\h^{i}$ exists, and that each $\h^{i}$ extends to a
continuous, piecewise polynomial function with respect to a finite
polyhedral decomposition of $A_{n-1}(X)_{\R}$.  Tchoudjem has given a similar combinatorial description of the cohomology of equivariant line bundles on regular compactifications of reductive groups \cite[Th\'eor\`eme 2.1]{Tchoudjem}; it should be interesting to investigate
whether the results and techniques of this paper extend to such varieties.

We apply our cohomology computations to give a formula for the
self-intersection number of a $T$-Cartier divisor.  For each 
$I \subset \Delta(1)$, let $\Delta_{I}$ be the fan
consisting of exactly those cones in $\Delta$ spanned by rays in $I$,
and let $\Delta_{I}(j)$ be the set of $j$-dimensional cones in
$\Delta_{I}$.  Define
\[
   \chi(\Delta_{I}) := \sum_{j=0}^{n} (-1)^{j} \cdot \# \Delta_{I}(j).
\]
In Section 2, we show that
\[
    \chi(\O(D)) = (-1)^{n} \sum_{P_{D,I} \mbox{ bounded}} 
    \chi(\Delta_{I}) \cdot (\# P_{D,I} \cap M).
\]
Using this formula for $\chi(\O(D))$ and asymptotic Riemann-Roch, we 
give a self-intersection formula for $T$-Cartier divisors.  When
$P_{D,I}$ is bounded, we write $\vol P_{D,I}$ for the volume of
$P_{D,I}$, normalized so that the smallest lattice simplex has unit
volume.

\begin{theorem}[Self-intersection formula] \label{self-intersection}
    Let $X$ be a complete $n$-dimensional toric variety and $D$ a
    $T$-Cartier divisor on $X$.  Then
    \[
	(D^{n}) = (-1)^{n} \cdot \sum_{P_{D,I} \mbox{ bounded}}
	 \chi(\Delta_{I}) \cdot \vol P_{D,I}.
    \]
\end{theorem}

\noindent When $X$ is smooth, Theorem \ref{self-intersection} is
closely related to a formula of Karshon and Tolman for the pushforward
of the top exterior power of a presymplectic form under the moment map
\cite{KT}.  In this case, the coefficient $(-1)^{n} \cdot
\chi(\Delta_{I})$ is equal to a winding number which gives the density
of the Duistermaat-Heckman measure on $P_{D,I}$.  

We conclude by proving an ``asymptotic converse'' to Serre vanishing
in the toric case.  From Serre vanishing we know that, for $D$ ample,
$h^{i} (mD) = 0$ for all $i > 0$ and $m \gg 0$.  The set of ample
divisors is open in $\Pic(X)_{\R}$, so the higher volume functions
vanish in a neighborhood of every ample divisor.  We prove the
converse for divisors on complete simplicial toric varieties.

\begin{theorem}[Asymptotic converse to Serre vanishing] \label{serre
converse}
    Let $D$ be a divisor on a complete simplicial toric variety.  Then
    $D$ is ample if and only if $\h^{i}$ vanishes identically in a
    neighborhood of $D$ in $\Pic (X)_{\R}$ for all $i > 0$.
\end{theorem}

\noindent The asymptotic converse to Serre vanishing does not hold in
general if $X$ is complete but not simplicial.  Fulton gives an
example of a complete, nonprojective toric threefold with no
nontrivial line bundles \cite[pp.  25-26, 72]{Fulton}.  For such an
$X$, $\Pic(X) = 0$ and all of the $\h^{i}$ vanish, but the zero
divisor is not ample.  We do not know whether the asymptotic converse
to Serre vanishing holds for nonsimplicial projective toric varieties.

On a toric variety, linear equivalence and numerical equivalence of
Cartier divisors coincide, so $\Pic(X)_{\R} = N^{1}(X)_{\R}$. 
Lazarsfeld asks whether, for a smooth complex projective variety $X$,
a divisor $D$ is ample if and only if the higher asymptotic
cohomological functions vanish in a neighborhood of the class of $D$
in $N^{1}(X)_{\R}$.\footnote{In the time since this article was written, de Fernex, Lazarsfeld, and the
second author have announced an affirmative answer to this question \cite[Theorem 4.7]{ELMNP}.}

\vspace{10 pt}

\noindent We thank R. Lazarsfeld, whose questions provided the 
starting point for this project, for his support and encouragement.

\section{Cohomology of $T$-Weil divisors}

By the cohomology groups of a Weil divisor $D$ on an algebraic variety
$X$, we always mean the sheaf cohomology groups $H^{i} (X, \O(D))$,
where $\O(D)$ is the sheaf whose sections over $U$ are the rational
functions $f$ such that $(\div f + D)|_{U}$ is effective.  When $X$ is
complete, we write $h^{i}(D)$ for the dimension of $H^{i}(X, \O(D))$.

In this section, for each $T$-Weil divisor $D$ on a toric variety, we
give a decomposition of the weight space $M_{\R}$ into finitely many
polyhedral regions such that the dimension of the $u$-graded piece of
the $i$-th cohomology group of $D$ is constant for all $u$ in each
region.  This decomposition can be deduced from local cohomology
computations in \cite[Theorem 2.7]{EMS}, but we present a proof using
different methods.  Our approach is a variation on the standard method
for computing the cohomology groups of $T$-Cartier divisors
\cite[Section 3.5]{Fulton}.

Let $X = X(\Delta)$ be an $n$-dimensional toric variety over a field
$k$, and let $\Delta(1)$ be the set of rays of $\Delta$.  Let $D =
\sum d_{\rho} D_{\rho}$ be a $T$-Weil divisor.  For each $I \subset
\Delta(1)$, define
\[ 
     P_{D,I} := \{ u \in M_{\R} : \< u, v_{\rho} \> \geq -d_{\rho}
     \mathrm{\ if\ and\ only\ if\ } \rho \in I \},
\] 
and let $\Delta_{I}$ be the subfan of $\Delta$ consisting of exactly
those cones whose rays are contained in $I$.  Note that
$P_{D,\Delta(1)}$ is the closed polyhedron usually denoted $P_{D}$,
each $P_{D,I}$ is a polyhedral region defined by an intersection of
halfspaces, some closed and some open, and $M_{\R}$ is their disjoint
union.  With $D$ fixed, for each $u \in M$, set
\[
    I_u := \{\rho \in {\Delta}(1) : \< u, v_{\rho} \> \geq -d_\rho \}.
\]

Recall that $H^{i}_{| \Delta_{I} |} (| \Delta |)$ denotes the
topological local cohomology group of $|\Delta|$ with support in
$|\Delta_{I}|$.  Here and throughout, all topological homology and
cohomology groups are taken with coefficients in $k$, the base field
of $X$.

\begin{proposition} \label{cohomology}
Let $X = X(\Delta)$ be a toric variety, $D$ a $T$-Weil divisor on 
$X$. Then 
\[
   H^{i} (X, \O(D)) \cong \bigoplus_{u \in M}  H^{i}_{| \Delta_{I_{u}} |} (| 
   \Delta |).
\]
\end{proposition}
 
\noindent \emph{Proof:} The \v{C}ech complex $C^{\bullet}$ that computes the
cohomology of $\O(D)$ is $M$-graded, and the $u$-graded piece is a 
direct sum of $u$-graded pieces of modules of sections of $\O(D)$ 
as follows:
\[
    C^{i}_{u} = \bigoplus_{\sigma_{0}, \ldots, \sigma_{i} \in \Delta}
    H^{0}( U_{\sigma_{0}} \cap \cdots \cap U_{\sigma_{i}},
    \O(D))_{u}.
\]     
Now $H^{0}( U_{\sigma_{0}} \cap \cdots \cap U_{\sigma_{i}},
\O(D))_{u}$ is isomorphic to $k$ if $\sigma_{0} \cap \cdots \cap
\sigma_{i}$ is in $\Delta_{I_{u}}$, and is zero otherwise.  In
particular,
\[
   H^{0}(U_{\sigma_{0}} \cap \cdots \cap U_{\sigma_{i}}, \O(D))_{u}
   \cong H^{0}_{|\Delta_{I_{u}}| \cap \sigma_{0} \cap \cdots \cap
   \sigma_{i}}(\sigma_{0} \cap \cdots \cap \sigma_{i}).
\]
A standard argument from topology \cite[Lemma p.75]{Fulton} shows that
the \v{C}ech complex $C^{\bullet}$ also computes $H^{i}_{|
\Delta_{I_{u}} |} (| \Delta |)$.  \hfill $\Box$

\begin{corollary} \label{polyhedral decomposition}
    Let $X = X(\Delta)$ be a toric variety, $D$ a $T$-Weil divisor on 
$X$. Then 
\[
   H^{i} (X, \O(D)) \cong \bigoplus_{I \subset \Delta(1)} \left( 
   \bigoplus_{u \in P_{D,I} \cap M} H^{i}_{|
   \Delta_{I} |} (| \Delta |) \right).
\]
\end{corollary}

\noindent \emph{Proof:} Since $P_{D,I} \cap M$ is exactly the set of
$u$ such that $I_{u} = I$, the corollary follows from Proposition
\ref{cohomology} by regrouping the summands.

\begin{proposition} \label{betti numbers}
   Let $X = X(\Delta)$ be a complete toric variety.  For $D$ a $T$-Weil 
   divisor on $X$,
    \[
	 h^{i}(D) = \sum_{P_{D,I} \mbox{ bounded} } h^{i}_{|
	 \Delta_{I} |} (N_{\R}) \cdot \# (P_{D,I} \cap M).
    \]
\end{proposition}

\noindent \emph{Proof:} When $X$ is complete, the support of $\Delta$
is all of $N_{\R}$, and $H^{i}(X, \O(D))$ is finite dimensional.  By
Corollary \ref{polyhedral decomposition},
$H^{i}_{|\Delta_{I}|}(|\Delta|)$ must vanish whenever $P_{D,I}$ is
unbounded, and the result follows.  \hfill $\Box$

\vspace{10 pt}

If $S$ is the unit sphere for some choice of coordinates on $N_{\R}$,
then $h^{i}_{|\Delta_{I}|}(N_{\R}) \cong \widetilde{h}_{n-i-1}(|\Delta_{I}|
\cap S)$ \cite[Exercise p.88]{Fulton}.  Therefore, Proposition
\ref{betti numbers} implies that computations of cohomology groups of
toric divisors can be reduced to computations of reduced homology
groups of finite polyhedral cell complexes and counting lattice points in
polytopes.  We will use the following lemma to show that the reduced
homology computations are not necessary if one is only interested in the
Euler characteristic $\chi(\O(D))$.  For any fan $\Sigma$, let
$\Sigma(j)$ denote the set of $j$-dimensional cones in $\Sigma$,
and define
\[
    \chi(\Sigma) := \sum_{j= 0}^{n} (-1)^{j} \cdot \# \Sigma(j).
\]    

\begin{lemma} \label{chi}
    Let $\Sigma$ be a fan in $N_{\R}$.  Then
    \[
	\sum_{i=0}^{n} (-1)^{i} \cdot h^{i}_{|\Sigma|}(N_{\R}) =
	(-1)^{n} \cdot \chi(\Sigma).
    \]
\end{lemma}

\noindent \emph{Proof:}  Let $S$ be the unit sphere for some choice of 
coordinates on $N_{\R}$.  Then
\begin{eqnarray*}
      \sum_{i=0}^{n} (-1)^{i} \cdot h^{i}_{|\Sigma|}(N_{\R}) & = & \sum_{i=0}^{n}
      (-1)^{i} \cdot \widetilde{h}_{n-i-1}(|\Sigma| \cap S) \\
          & = & (-1)^{n} + \sum_{i = 0}^{n} (-1)^{i} \cdot 
          h_{n-i-1}(|\Sigma| \cap S).
\end{eqnarray*}
Setting $j = n - i$, and then using the correspondence between the 
$j-1$-dimensional cells in $|\Sigma| \cap S$ and the $j$-dimensional 
cones in $\Sigma$, we have
\begin{eqnarray*}
          (-1)^{n} + \sum_{i = 0}^{n} (-1)^{i} \cdot 
          h_{n-i-1}(|\Sigma| \cap S) & = & (-1)^{n} + \sum_{j=0}^{n} 
          (-1)^{n-j} \cdot h_{j-1}(|\Sigma| \cap S). \\
	  & = & (-1)^{n} \cdot \sum_{j=0}^{n} (-1)^{j} \cdot \# \Sigma(j).  
\end{eqnarray*}

\vspace{-25 pt} \hfill $\Box$

\vspace{25 pt}

\begin{proposition} \label{euler characteristic}
    Let $D$ be a $T$-Weil divisor on a complete $n$-dimensional toric
    variety.  Then
    \[
	\chi(\O(D)) = (-1)^{n} \cdot \sum_{P_{D,I} \mbox{ bounded}} 
	\chi(\Delta_{I}) \cdot \# (P_{D,I} \cap M).
    \]	
\end{proposition}

\noindent \emph{Proof:} The proposition follows immediately from Proposition
\ref{betti numbers} and Lemma~\ref{chi}.  \nolinebreak \hfill $\Box$

\section{Asymptotic cohomological functions and the self-intersection 
formula}

\begin{definition} Let $X$ be a complete
$n$-dimensional toric variety.  The $i$-th asymptotic cohomological
function $\h^{i} : A_{n-1}(X) \rightarrow \R$ is defined by 
    \[
       \h^{i}(D) = \lim_{m} \frac {h^{i}(mD)} {m^{n} / n!}.
    \]    
\end{definition}

For a bounded polyhedral region $P \subset M_{\R}$, let $\vol P$ 
denote the volume of $P$, normalized so that the smallest lattice 
simplex has unit volume.  Note that
\[
    \vol P = \lim_{m} \frac{ \# mP \cap M } {m^{n} / n!}.
\]    

\begin{proposition} \label{higher volumes}
    Let $D$ be a $T$-Weil divisor on a complete toric variety $X = 
    X(\Delta)$.  Then
    \[
       \h^{i}(D) = \sum_{P_{D,I} \mbox{ bounded}} h^{i}_{| \Delta_{I}
       |} (N_{\R}) \cdot \vol P_{D,I},
    \]
\end{proposition}

\noindent \emph{Proof:} For all $I \subset \Delta(1)$, and for all positive 
integers $m$, $P_{mD,I} = m P_{D,I}$.  The proposition therefore 
follows immediately from the definition of $\h^{i}$ and Proposition 
\ref{betti numbers}.  \hfill $\Box$

\begin{corollary} Let $X$ be a complete
$n$-dimensional toric variety.  Then $\h^{i}$ extends to a
continuous, piecewise polynomial function with respect to a finite
polyhedral decomposition of $A_{n-1}(X)_{\R}$. 
\end{corollary}

\noindent \emph{Proof:} The set of $I$ such that $P_{D,I}$ is bounded does
not depend on $D$.  Indeed, $P_{D,I}$ is bounded if and only if there
is no hyperplane in $N_{\R}$ separating the rays in $I$ from the rays
in $\Delta(1) \smallsetminus I$.  The result then follows from Proposition
\ref{higher volumes} since, for each such $I$, $\vol P_{D,I}$ extends
to a continuous, piecewise polynomial function with respect to a
finite polyhedral decomposition of $A_{n-1}(X)_{\R}$.  \hfill $\Box$

\vspace{10 pt}

\noindent \textbf{Theorem 1 (Self-intersection formula)} \emph{
    Let $X$ be a complete $n$-dimensional toric variety and $D$ a
    $T$-Cartier divisor on $X$.  Then }
    \[
	(D^{n}) = (-1)^{n} \cdot \sum_{P_{D,I} \mbox{ bounded}}
	 \chi(\Delta_{I}) \cdot \vol P_{D,I}.
    \]
    
\vspace{5 pt}    

\noindent \emph{Proof:} By asymptotic Riemann-Roch
\cite[VI.2]{Kollar}, when $D$ is Cartier,
\[
     (D^{n}) = \lim_{m} \frac{\chi(\O(mD))}{m^{n}/n!}.
\]
The theorem then follows from Proposition \ref{euler
characteristic}. \hfill $\Box$

\section{Asymptotic converse to Serre vanishing}

We begin by briefly recalling the Gelfand-Kapranov-Zelevinsky (GKZ)
decomposition introduced by Oda and Park \cite{OP} and a few of its
basic properties.  Assume that $X$ is complete.  The GKZ decomposition
is a fan whose support is the effective cone in $A_{n-1}(X)_{\R}$ and
whose maximal cones are in 1-1 correspondence with the simplicial fans
$\Sigma$ in $N_{\R}$ such that $\Sigma(1) \subset \Delta(1)$ and
$X(\Sigma)$ is projective.  We call the interior of a maximal GKZ cone
a \emph{GKZ chamber}, and write $\gamma_{\Sigma}$ for the GKZ chamber
corresponding to $\Sigma$.  If $D$ is a $T$-Weil divisor whose class
$[D]$ lies in $\gamma_{\Sigma}$, then $\Sigma$ is the normal fan to
$P_{D}$.  This property fully characterizes the GKZ decomposition.  We
will need the following basic property relating divisors in
$\gamma_{\Sigma}$ to divisors on $X(\Sigma)$: if $f$ denotes the
birational map from $X$ to $X(\Sigma)$ induced by the identity on $N$,
then the birational transform $f_{*}(D)$ is ample on $X(\Sigma)$, and
$P_{f_{*}(D)} = P_{D}$.  See the appendix for proofs and for a more
detailed discussion of the GKZ decomposition in the language of toric
divisors.

\begin{lemma} \label{partials} Let $\gamma_{\Sigma}$ be a GKZ chamber,
and let $f$ be the birational map from $X = X(\Delta)$ to $X(\Sigma)$
induced by the identity on $N$.  Let $D_{1}, \ldots, D_{r}$ be
distinct prime $T$-invariant divisors on $X$ corresponding to rays
$\rho_{1}, \ldots, \rho_{r} \in \Delta$, respectively.  For $D$ a
$T$-Weil divisor with $[D] \in \gamma_{\Sigma}$,
\[
    \frac{ \partial^{r} \h^{0}} {\partial D_{1} \cdots \partial 
    D_{r}} (D) = \frac{n!}{(n-r)!} \cdot \left( f_{*}(D)^{n-r} \cdot f_{*}(D_{1}) \cdot 
    \ldots \cdot f_{*}(D_{r}) \right).
\]
In particular, $\frac{ \partial^{r} \h^{0}} {\partial D_{1} \cdots
\partial D_{r}}$ is strictly positive on $\gamma_{\Sigma}$ if
$\rho_{1}, \ldots, \rho_{r}$ span a cone in $\Sigma$ and vanishes 
identically on $\gamma_{\Sigma}$ otherwise.
\end{lemma}

\noindent \emph{Proof:} Suppose $r = 1$.  Since $f_{*}(D)$ is ample
and $P_{f_{*}(D)} = P_{D}$ for $D$ in $\gamma_{\Sigma}$, $\h^{0}$ is given
on $\gamma_{\Sigma}$ by $D \mapsto (f_{*}(D)^{n})$.  Therefore,
\begin{eqnarray*}
    \frac {\partial \h^{0}} {\partial D_{1}} & = & \lim_{\epsilon
    \rightarrow 0} \frac {\left( f_{*}(D + \epsilon D_{1})^{n} \right)
    - \left( f_{*}(D)^{n} \right)} {\epsilon}.  \\
         & = & n \left( f_{*}(D)^{n-1} \cdot f_{*}(D_{1}) \right).
\end{eqnarray*}
The general case follows by a similar computation and induction on
$r$.  The last statement follows from the formula, since $f_{*}(D)$ is
ample and $f_{*}(D_{1}) \cdot \ldots \cdot f_{*}(D_{r})$ is an
effective cycle if $\rho_{1}, \ldots, \rho_{r}$ span a cone in
$\Sigma$ and is zero otherwise \cite[Chapter 5]{Fulton}.  \hfill
$\Box$

\vspace{10 pt}

\noindent \textbf{Theorem 2 (An asymptotic converse to Serre 
vanishing)} \emph{ 
    Let $D$ be a divisor on a complete simplicial toric variety.  Then
    $D$ is ample if and only if $\h^{i}$ vanishes identically in a
    neighborhood of $D$ in $\Pic (X)_{\R}$ for all $i > 0$.}

\vspace{ 10 pt }    

\noindent \emph{Proof:} Since the limits in the definition of the $\h^{i}$
exist, by asymptotic Riemann-Roch, for $D$ a \Q -Cartier divisor,
$(D^{n}) = \sum_{i=0}^{n} (-1)^{i} \cdot \h^{i}(D).$ Therefore, if
$\h^{i}$ vanishes in a neighborhood of $D$ for all $i > 0$, then
$\h^{0}$ agrees with the top intersection form in a neighborhood of
$D$.  To prove Theorem \ref{serre converse}, we will prove the
stronger fact that if $\h^{0}$ agrees with the top intersection form
in a neighborhood of $D$, then $D$ is ample.  It will suffice to show
that if $\gamma_{\Sigma}$ is a GKZ chamber and $\h^{0}(D) = (D^{n})$ for
$[D] \in \gamma_{\Sigma}$, then $\Sigma = \Delta$.

Suppose $\gamma_{\Sigma}$ is a GKZ chamber and $\h^{0}(D) = (D^{n})$ for
$[D] \in \gamma_{\Sigma}$.  Let $\rho_{1}, \ldots, \rho_{n}$ be rays
spanning a maximal cone $\sigma \in \Delta$.  It will suffice to show
that $\rho_{1}, \ldots, \rho_{n}$ span a cone in $\Sigma$.  On
$\gamma_{\Sigma}$, since $\h^{0}$ agrees with the top intersection
form,
\begin{eqnarray*}
    \frac{ \partial^{n} \h^{0}} {\partial D_{1} \cdots \partial
    D_{n}} & = & n! \cdot \left( D_{1} \cdot \ldots \cdot D_{n} \right) \\
                 & = & n! \cdot \mult (\sigma)^{-1}.
\end{eqnarray*}
In particular, $\frac{ \partial^{n} \h^{0}} {\partial D_{1} \cdots
\partial D_{n}}$ does not vanish identically on $\gamma_{\Sigma}$.  By
Lemma \ref{partials}, $\rho_{1},$ \ldots, $\rho_{n}$ span a cone in
$\Sigma$, as required.  \hfill $\Box$

\section{Appendix: Gelfand-Kapranov-Zelevinsky Decompositions}

In this appendix we give a self-contained account of the GKZ
decompositions of Oda and Park \cite{OP}, in the language of toric
divisors.

A possibly degenerate fan in $N$ is a finite collection $\Sigma$ of
convex (not necessarily strongly convex) rational polyhedral cones in
$N_{\R}$ such that every face of a cone in $\Sigma$ is in $\Sigma$,
and the intersection of any two cones in $\Sigma$ is a face of each. 
The intersection of all of the cones in $\Sigma$ is the unique linear
subspace $L_{\Sigma} \subset N_{\R}$ that is a face of every cone in
$\Sigma$; we say that $\Sigma$ is degenerate if $L_{\Sigma}$ is not
zero.  Associated to $\Sigma$ is a toric variety $X_{\Sigma}$ of
dimension $\dim N_{\R} - \dim L_{\Sigma}$, whose torus is
$T_{N/L_{\Sigma} \cap N}$, and the $T$-Cartier divisors on
$X_{\Sigma}$ correspond naturally and bijectively to the piecewise
linear functions on $|\Sigma|$ whose restriction to $L_{\Sigma}$ is
identically zero.

Let $X = X(\Delta)$ be an $n$-dimensional toric variety and assume,
for simplicity, that $|\Delta |$ is convex and $n$-dimensional.  Let
$D = \sum d_{\rho} D_{\rho}$ be an effective $T$-\Q -Weil divisor and let
$P_D = \{ u \in M : \<u, v_{\rho}\>\geq -d_{\rho}\}$ be the
polyhedron associated to $D$.  From $P_D$ one constructs the possibly
degenerate normal fan $\Sigma _D$, whose support is $|\Delta|$ and
whose cones are in one to one order reversing correspondence with the
faces of $P_D$; the cone corresponding to a face $Q$ is
\[
    \sigma _Q = \{ v\in |\Delta| : \< u, v\> \geq \< u', v\> \textrm{
    for all } u\in P_D \textrm{ and } u'\in Q\}.
\]
Note that $\sigma _Q$ is positively spanned by those rays $\rho \in
\Delta (1)$ such that the affine hyperplane $\< u, v_{\rho} \> =
-d_{\rho}$ contains $Q$.

We define a convex piecewise linear function $\Xi_{D}$ on $|\Delta|$
by
\[
    \Xi_{D}(v) = \min \{ \< u, v \> : u \in P_{D} \}.
\]
The maximal cones of $\Sigma _D$ are the maximal domains of linearity
of $\Xi_D$.  When $D$ is \Q -Cartier and ample, $\Sigma_{D} = \Delta$
and $\Xi_{D} = \Psi_{D}$ is the piecewise linear function usually
associated to $D$ \cite[Section 3.3]{Fulton}.  It follows from the
definition of $\Xi_D$ that
\begin{equation}\label{restriction to rays}
 \Xi_{D}(v_{\rho}) \geq -d_{\rho},
\end{equation}
with equality for those $\rho$ such that the affine hyperplane $\< u,
v_{\rho} \> = -d_{\rho}$ contains a face of $P_{D}$.  Let $ I_{D}
\subset \Delta(1)$ be the set of rays for which equality does not hold
in (\ref{restriction to rays}).

\begin{definition}[GKZ cones] \label{GKZ cones} Let $\Sigma$ be a
possibly degenerate fan whose support is $|\Delta|$, such that
$X_{\Sigma}$ is quasiprojective, and for which there is a set of rays
$I \subset \Delta(1)$ such that every cone in $\Sigma$ is positively
spanned by rays in $\Delta(1) \smallsetminus I$.  The GKZ cone
$\gamma_{\Sigma,I}$ is defined to be
 \[
    \gamma_{\Sigma,I} := \{ [D] \in A_{n-1}(X)_{\Q} : \Sigma \textrm{
    refines } \Sigma _{D} \textrm{ and } I_D\subseteq I\}.
 \]
\end{definition}   

\vspace{3 pt}

\noindent The GKZ cone $\gamma_{\Sigma,I}$ is well-defined since $\Sigma_{D}$
and $I_{D}$ depend only on the linear equivalence class of $D$.

\vspace{10 pt}

\noindent \textbf{GKZ Decomposition Theorem} \cite[Theorem 3.5]{OP}
\emph{ The GKZ cone $\gamma_{\Sigma,I}$ is a rational polyhedral cone
of dimension $ \dim \Pic (X_{\Sigma})_{\Q} + \# I$.  The set of GKZ
cones is a fan whose support is the effective cone in
$A_{n-1}(X)_{\R}$, and the faces of $\gamma_{\Sigma, I}$ are exactly
those $\gamma_{\Sigma',I'}$ such that $\Sigma$ refines $\Sigma'$ and
$I' \subset I$.  }

\vspace{10 pt}

It follows from the theorem that the maximal GKZ cones are in 1-1
correspondence with the nondegenerate simplicial fans $\Sigma$ in
$N_{\R}$ such that $\Sigma(1) \subset \Delta(1)$, $|\Sigma| =
|\Delta|$, and $X(\Sigma)$ is quasiprojective.  Indeed, if $[D]$ is in
$\gamma_{\Sigma,I}$ then $\dim P_{D} = n - \dim L_{\Sigma}$.  In
particular, if $\Sigma$ is degenerate then $\gamma_{\Sigma,I}$ is in
the boundary of the effective cone and $\h^{0}|_{\gamma_{\Sigma,I}}$
is identically zero.  If $\Sigma$ is nondegenerate, then
\[
   \dim \gamma_{\Sigma,I} \leq \# \Sigma(1) -n + \# I \leq \#\Delta(1)
   -n,
\]   
with equalities if and only if $\Sigma$ is simplicial and $I =
\Delta(1) \smallsetminus \Sigma(1)$.  The interiors of the maximal
cones of the GKZ decomposition are called GKZ chambers, and we
write $\gamma_{\Sigma}$ for the GKZ chamber corresponding to $\Sigma$.

In order to prove the GKZ Decomposition Theorem, we need a few basic
tools relating divisors on $X$ to divisors on $X_{\Sigma}$, where
$\Sigma$ is a possibly degenerate fan in $N_{\R}$ whose support is
$|\Delta|$.  Let $\phi_{\Sigma}$ be the map taking a $T$-\Q -Cartier
divisor $D$ on $X_{\Sigma}$ to the $T$-\Q -Weil divisor
$\phi_{\Sigma}(D)$ on $X$, where
\[
    \phi_{\Sigma}(D) = \sum_{\rho \in \Delta(1)} -\Psi_{D}(v_{\rho}) 
    D_{\rho}.
\]
Note that $\phi_{\Sigma}$ respects linear equivalence and induces an
injection of $\Pic(X_{\Sigma})$ into $A_{n-1}(X)$.  The map
$\phi_{\Sigma}$ may be realized geometrically as follows.  Let
$\widetilde{X}$ be the toric variety corresponding to the smallest
common refinement of $\Delta$ and $\Sigma$.  The identity on $N$
induces morphisms $p_{1}$ and $p_{2}$ from $\widetilde{X}$ to $X$ and
to $X_{\Sigma}$, respectively.  Then $\phi_{\Sigma} = p_{1_{*}} \circ
p_{2}^{*}$.  

The following lemma will be used to prove the GKZ Decomposition
Theorem.  It also shows that $\gamma _{\Sigma, I}$ is equal to the
cone $\mathrm{cpl}(\Sigma, \Delta(1)\smallsetminus I)$ defined in
\cite[Section~3]{OP}.

\begin{lemma}\label{GKZ lemma} Let $D = \sum d_{\rho}D_{\rho}$ be a
$T$-\Q -Weil divisor.  The following are equivalent:
  \begin{enumerate} \renewcommand{\theenumi}{\roman{enumi}}
  \item  The GKZ cone $\gamma_{\Sigma,I}$ contains $[D]$.  
  
  \item  There is a convex function $\Xi$ that is linear on each maximal
  cone of $\Sigma$ such that $\Xi (v_{\rho}) \geq -d_{\rho}$ for all
  $\rho \in \Delta(1)$, with equality when $\rho \notin I$. 
  
  \item There is a divisor $\widetilde{D}$ linearly equivalent to $D$ and
  a decomposition $ \widetilde{D} = \phi_{\Sigma}(D') + E$ such that $D'$
  is a nef $T$-\Q -Cartier divisor on $X_{\Sigma}$, $P_{\widetilde{D}} =
  P_{D'}$, and $E$ is an effective divisor whose support is contained
  in $\bigcup_{\rho \in I} D_{\rho}$.
  \end{enumerate}
\end{lemma}

\noindent \emph{Proof:} If $[D]$ is in $\gamma_{\Sigma,I}$, then (ii)
holds for $\Xi = \Xi_{D}$.  If (ii) holds, then choose $u \in M_{\Q}$
such that $\Xi|_{L_{\Sigma}} = u|_{L_{\Sigma}}$.  Let $\widetilde{D} =
D + \sum_{\rho} \< u, v_{\rho} \> D_{\rho}$, and let $D'$ be the \Q
-Cartier divisor on $X_{\Sigma}$ corresponding to $\Xi - u$.  Since
$\Xi$ is convex, $D'$ is nef, and (iii) holds with $E = \sum_{\rho}
\left( \Xi(v_{\rho}) + d_{\rho} \right) D_{\rho}$.  Also, since
$\Psi_{D'} = \Xi_{D} - u = \Xi_{\widetilde{D}}$, we have $P_{D'} =
P_{\widetilde{D}}$.

It remains to show that (iii) implies (i).  Replacing $D$ by
$\widetilde{D}$ if necessary, we may assume $D = \phi_{\Sigma}(D') + E$,
where $D'$ is a nef \Q -Cartier divisor on $X_{\Sigma}$ and $E$ is an
effective divisor whose support is contained in $\bigcup_{\rho \in I}
D_{\rho}$.  We must show that $\Sigma$ refines $\Sigma_{D}$ and
$\Xi_{D}(v_{\rho}) = -d_{\rho}$ for $\rho \not \in I$. Since $P_{D} =
P_{D'}$, and since $D'$ is nef, $\Xi_{D} = \Psi_{D'}$,
which is linear on each cone of $\Sigma$.  Hence $\Sigma_{D}$ is
refined by $\Sigma$.  Since the support of $E$ is contained in
$\bigcup_{\rho \in I} D_{\rho}$, we also have $\Xi_{D}(v_{\rho}) =
-d_{\rho}$ for $\rho \not \in I$, as required.  \hfill $\Box$ 

\begin{corollary} \label{relative interiors}
  Let $D$ be a $T$-Weil divisor on $X$.  Then $[D]$ is in the relative
  interior of $\gamma_{\Sigma,I}$ if and only if $\Sigma_{D} = \Sigma$
  and $I_{D} = I$.
\end{corollary} 

\noindent \emph{Proof:} The decomposition in Lemma \ref{GKZ lemma} part (iii)
is essentially unique; if $ \widetilde{D}$ is replaced by $\widetilde{D} +
\sum \<u, v_{\rho}\> D_{\rho}$, where $u|_{L_{\Sigma}} = 0$, then $D'$
is replaced by the divisor corresponding to $\Psi_{D'} - u$ and $E =
\sum_{\rho \in I} e_{\rho} D_{\rho}$ remains fixed.  It follows that
the map $[D] \mapsto (D', (e_{\rho})_{\rho \in I})$ gives an
isomorphism
\[
   \gamma_{\Sigma,I} \stackrel{\sim}{\longrightarrow} 
   \mathrm{Nef}(X_{\Sigma}) \times \R^{I}_{\geq 0}.
\]
Taking relative interiors gives $\gamma_{\Sigma,I}^{\circ}
\stackrel{\sim}{\longrightarrow} \mathrm{Ample}(X_{\Sigma}) \times
\R^{I}_{> 0}.$  Therefore $[D]$ is in the relative interior of 
$\gamma_{\Sigma,I}$ if and only if $\Xi_{D}$ is strictly convex with 
respect to $\Sigma$ and the inequality in (\ref{restriction to rays}) 
is strict exactly when $\rho \in I$. \hfill $\Box$

\begin{corollary} \label{distinct volumes}
  Suppose $X$ is complete, and let $\gamma_{\Sigma,I}$ be a GKZ cone,
  with $\Sigma$ nondegenerate.  Let $f$ be the birational map from $X$
  to $X_{\Sigma}$ induced by the identity on $N$.  If $[D] \in
  \gamma_{\Sigma,I}$ (resp.  $[D] \in \gamma_{\Sigma,I}^{\circ}$),
  then $P_{f_{*}(D)} = P_{D}$ and $f_{*}(D)$ is nef (resp.  $f_{*}(D)$
  is ample).  In particular, $\h^{0}|_{ \gamma_{\Sigma,I} }$ is given
  by $[D] \mapsto (f_{*}(D)^{n})$.
\end{corollary}   

\noindent \emph{Proof:} Since $\Sigma$ is nondegenerate, we can take
$\widetilde{D} = D$ and let $D = \phi_{\Sigma}(D') + E$ be the
decomposition in Lemma \ref{GKZ lemma} part (iii).  Then, since $P_{D}
= P_{D'},$ $\h^{0}(D) = \h^{0}(D')$.   Furthermore, since $D'$ is nef on
$X_{\Sigma}$ (and is ample if $[D] \in \gamma_{\Sigma,I}^{\circ})$,
$\h^{0}(D') = ((D')^{n})$.  We claim that $D' = f_{*}(D)$.  
Indeed,
\[  
   D' = \sum_{\rho \in \Sigma(1)} - \Xi_{D}(v_{\rho}) f_{*}(D_{\rho})
   = f_{*}(D),
\]
and the result follows.  \hfill $\Box$

\vspace{10 pt}

\begin{corollary}
    The volume function of a complete toric variety is given by distinct 
    polynomials on distinct GKZ chambers.
\end{corollary}

\noindent \emph{Proof:} Let $\gamma_{\Sigma}, \gamma_{\Sigma'}$ be 
distinct GKZ chambers.  Let $\rho_{1}, \ldots, \rho_{n}$ be rays 
spanning a maximal cone in $\Sigma$ that is not in $\Sigma'$.  
By Lemma \ref{partials}, $\frac {\partial^{n} \h^{0} } {\partial 
D_{1} \cdots \partial D_{n}}$ vanishes identically on 
$\gamma_{\Sigma'}$, but not on $\gamma_{\Sigma}$.

\vspace{10 pt}

\noindent \emph{Proof of GKZ Decomposition Theorem:} First, we claim
that $[D]$ is in $\gamma_{\Sigma,I}$ if and only if for each maximal
cone $\sigma \in \Sigma$, for each collection of linearly independent
rays $\rho_{1}, \ldots, \rho_{n}$ in $\Delta(1) \smallsetminus I$ that
are contained in $\sigma$, and for each $\rho \in \Delta(1)$ with
$v_{\rho} = a_{1} v_{\rho_{1}} + \cdots + a_{n} v_{\rho_{n}}$, we have
\begin{equation} \label{GKZ inequalities}
    - a_{1}d_{\rho_{1}} - \cdots - a_{n}d_{\rho_{n}} \left\{
    \begin{array}{ll} = & -d_{\rho} \mbox{ if } \rho \subset \sigma
    \mbox{ and } \rho \not \in I. \\ \geq & -d_{\rho} \mbox
    { otherwise.} \end{array} \right.
\end{equation}
There are only finitely many such conditions, and all of the
coefficients are rational, so the claim implies that each
$\gamma_{\Sigma,I}$ is a convex rational polyhedral cone.  Suppose
(\ref{GKZ inequalities}) holds.  Let $u_{\sigma} \in M_{\Q}$ be such
that $\< u_{\sigma}, v_{\rho_{i}} \> = -d_{\rho_{i}}$ for $1 \leq i
\leq n$.  The equalities in (\ref{GKZ inequalities}) ensure that the
$u_{\sigma}$ glue together to give a continuous piecewise linear
function $\Xi$ on $|\Sigma|$, where $\Xi|_{\sigma} = u_{\sigma}$, such
that $\Xi(v_{\rho}) = -d_{\rho}$ for $\rho \not \in I$.  The
inequalities in (\ref{GKZ inequalities}) guarantee that $\Xi$ is
convex.  By part (ii) of Lemma \ref{GKZ lemma}, it follows that $[D]$
is in $\gamma_{\Sigma,I}$.  Conversely, if $[D]$ is in
$\gamma_{\Sigma, I}$, then we have a $\Xi$ as in part (ii) of Lemma
\ref{GKZ lemma}.  Say $\Xi|_{\sigma} = u_{\sigma}$.  Then the left
hand side of (\ref{GKZ inequalities}) is equal to $\< u_{\sigma},
v_{\rho} \>$ and the desired equalities and inequalities follow from 
the choice of $\Xi$.

It follows from Corollary \ref{relative interiors} that the effective
cone in $A_{n-1}(X)_{\R}$ is the disjoint union of the relative
interiors of the GKZ cones, and that $\dim \gamma_{\Sigma,I} = \dim
\Pic(X_{\Sigma})_{\Q} + \#I.$ Any finite collection of rational
polyhedral cones such that every face of a cone in the collection is
in the collection, and such that the relative interiors of the cones
are disjoint, is a fan.  The faces of a cone in such a collection are
precisely the cones in the collection that it contains.  Therefore, to
prove the theorem, it remains only to show that every face of a GKZ
cone is a GKZ cone.  Let $\gamma_{\Sigma,I}$ be a GKZ cone, let
$\rho_{1}, \ldots, \rho_{n}$ be linearly independent rays contained in
a maximal cone $\sigma \in \Sigma$, let $\rho \in \Delta(1)$ with
$v_{\rho} = a_{1} v_{\rho_{1}} + \cdots + a_{n}v_{\rho_{n}}$, and let
$\tau \preceq \gamma_{\Sigma,I}$ be the face where equality holds in
(\ref{GKZ inequalities}).  If $\rho \subset \sigma$, then $\tau =
\gamma_{\Sigma, I \smallsetminus \{ \rho \} }$.  If $\rho \not \subset
\sigma$ then consider the set of convex cones $\sigma'$ in $\N_{\R}$
which are unions of maximal cones in $\Sigma$, which contain $\sigma$
and $\rho$, and are such that $X(\Sigma')$ is quasiprojective, where
$\Sigma'$ is the fan whose maximal cones are $\sigma'$ and all of the
maximal cones of $\Sigma$ that are not contained in $\sigma'$.  This
set is nonempty since it contains $|\Sigma|$, and since it is closed
under intersections it must contain a minimal element
$\overline{\sigma}$.  Let $\overline{\Sigma}$ be the corresponding
fan.  Then $\tau = \gamma_{\overline{\Sigma}, I}$.  \nolinebreak
\hfill $\Box$

\end{document}